\documentclass[11pt,a4paper]{amsart}
\usepackage{pagesize}
\usepackage{amssymb}
\usepackage{amsmath}
\theoremstyle{plain}
\newtheorem{Thm}{Theorem}
\newtheorem{Cor}{Corollary}
\newtheorem{Lem}{Lemma}
\newtheorem{Prop}{Proposition}

\theoremstyle{definition}
\newtheorem{Def}{Definition}
\theoremstyle{remark}

\numberwithin{equation}{section}

\pagesize{verbose=,paper=a4,width=146.67truemm,height=248.2truemm}
\begin{document}
\title{The pluripolar hull of a graph and fine analytic continuation}
\author{Tomas Edlund and Burglind J\"oricke}
\address{Tomas Edlund \\Department of Mathematics \\ Uppsala University\\Box 480\\SE-751 06 Uppsala \\ Sweden}
\email{tomas@math.uu.se}
\address{ Burglind J\"oricke \\Department of Mathematics \\ Uppsala University\\Box 480\\SE-751 06 Uppsala \\ Sweden}
\email{joericke@math.uu.se}
\date{\today}

\begin{abstract}
We show that if the graph of a bounded analytic function in the unit disk $\mathbb D$ is not complete pluripolar in $\mathbb C^2$ then the projection of the closure of its pluripolar hull contains a fine neighborhood of a point $p \in \partial \mathbb D$. On the other hand we show that if an analytic function $f$ in $\mathbb D$ extends to a function $\mathcal{F}$ which is defined on a fine neighborhood of a point $p \in \partial \mathbb D$ and is finely analytic at $p$ then the pluripolar hull of the graph of $f$ contains the graph of $\mathcal{F}$ over a smaller fine neighborhood of $p$. We give several examples of functions with this property of fine analytic continuation. 
As a corollary we obtain new classes of analytic functions in the disk which have non-trivial pluripolar hulls, among them $C^\infty$ functions on the closed unit disk which are nowhere analytically extendible and have infinitely-sheeted pluripolar hulls. Previous examples of functions with non-trivial pluripolar hull of the graph have fine analytic continuation.
\end{abstract}
\maketitle

\section{Introduction}\label{intro}\noindent
A subset $E$ of a domain $\Omega \subset \mathbb C^N$ is called pluripolar in $\Omega$, if for all $z \in E$ there exist a connected neighborhood $U_z$ of $z$ in $\Omega$ and a plurisubharmonic function $u(z,w) \not\equiv  -\infty$ defined on $U_z$ such that
\begin{displaymath}
E \cap U_z \subset \{(z,w)\in U_z :  u(z,w)= -\infty\}.
\end{displaymath}
By Josefson's theorem  (see \cite{josefson}), a set $E \subset \mathbb C^N$ is pluripolar if and only if there exists a globally defined plurisubharmonic function $u(z,w)$ such that
\begin{displaymath}
E \subset \{(z,w)\in \mathbb C^N : u(z,w)= -\infty\}.
\end{displaymath}
In other words a pluripolar set is a subset of the $-\infty$-locus of a globally defined plurisubharmonic function. Pluripolar sets are the exceptional sets in pluripotential theory. This motivates the interest in understanding the structure of pluripolar sets. A set $E \subset \Omega$ is called complete pluripolar in $\Omega$ if $E$ is the exact $-\infty$-locus of a plurisubharmonic function defined in $\Omega$. On the contrary, some subsets $E \subset \Omega$ (e.g. proper open subsets of connected analytic submanifolds) have the property that any plurisubharmonic function which is $-\infty$ on $E$ is $-\infty$ on a larger set. This leads to the notion of the {\em{pluripolar hull}} $E_{\Omega}^*$ (see \cite{hull}) of a pluripolar subset $E \subset \Omega$, 
\begin{displaymath}
E_{\Omega}^* \overset{def}{=} \bigcap \{z \in \Omega : u(z)= -\infty \},
\end{displaymath}
where the intersection is taken over {\em{all }} plurisubharmonic
functions in $\Omega$ which equal $-\infty$ on $E$.
In general, it is difficult to describe the pluripolar
hull of a given set $E$. Initiated by a paper of Saddulaev (\cite{Sad}) the pluripolar hull of graphs of certain analytic functions has been studied in a number of papers (see e.g. \cite{Ediwie},\cite{Ediwie2},\cite{Ediwie3}, \cite{hull}, \cite{js1},\cite{js2},\cite{w1} and \cite{w2}).

For a subset $A$ of the complex plane $\mathbb C$ and a complex valued function $f$ on $A$ we denote by $\Gamma_f(A)$ the graph of $f$ over $A$,
\begin{equation*}
 \Gamma_f(A) = \{ (z,w) \in \mathbb C^2 :z \in A, w =f(z) \}.
\end{equation*}
Let $f$ be a holomorphic function in the unit disk $\mathbb D \subset \mathbb C$. Clearly $\Gamma_f(\mathbb D)$ is a pluripolar set. It is a natural attempt to relate non-triviality of the pluripolar hull of $\Gamma_f(\mathbb D)$ to the existence of analytic continuation or various kinds of generalized analytic continuation of $f$ across some part of $\partial \mathbb D$. In \cite{analdisc} Levenberg, Martin and Poletsky conjectured that if
 $f$ is analytic in $\mathbb D$ and $f$ does not extend
 holomorphically across $\partial \mathbb D$, then $\Gamma_f(\mathbb D)$ is
 complete pluripolar. This conjecture was disproved in
 \cite{Ediwie2}. In \cite{js1}, Siciak noticed that the function in
 \cite{Ediwie2} possesses pseudocontinuation across a subset of the circle of
 full measure and showed that the pluripolar hull of its graph
 contains the graph of the pseudocontinuation. He also noticed that if an analytic function $f$ in $\mathbb D$ admits
pseudocontinuation through a set $E$ of positive measure on the circle
and the graph of the non-tangential limits is in the pluripolar hull
of  $\Gamma_f(\mathbb D)$ then also the graph of the pseudocontinuation is in the
mentioned pluripolar hull. In \cite{js1} he showed by an example that the existence of pseudocontinuation of the function $f$ is not necessary for non-triviality of the pluripolar hull of $\Gamma_f(\mathbb D)$.

The notion of fine analytic continuation seems to us better adapted to understand pluripolar completeness of graphs.

Recall that the fine topology was introduced by Cartan (see e.g. \cite{Bre}) as the weakest topology for which all subharmonic functions are continuous. A neighborhood basis of a point in this topology consists of sets which differ from a Euclidean neighborhood of this point by a set which is thin at this point. Thin sets were introduced by Brelot. A set $F \subset \mathbb C$ is thin at a point $\xi$, if either $\xi$ is not in its closure $\overline{F}$ or $\xi \in \overline{F}$ and there exists a subharmonic function $\mathcal{V}$ in a neighborhood of $\xi$ such that $\overline{\lim}_{z \in F, z \rightarrow \xi} \mathcal{V}(z) < \mathcal{V}(\xi)$. One can always choose $\mathcal{V}$ in such a way that the limit on the left equals $-\infty$. 
For a point $p \in \mathbb C$ and a positive number $r$ we denote by $D(p,r)$ the open disk of radius $r$ and center $p$.

By a closed fine neighborhood $V$ of $p$ we mean a connected closed set which has the form $B \smallsetminus U$ for some connected closed set $B$ and an open set $U \subset \mathbb C$ which is thin at $p$. Note that $U$ can be taken simply connected. We will often consider $B = \overline{D(p,r)}$ for some $r > 0$.  
\begin{Def}\label{Def2}
A continuous function $\mathcal{F}$ on a closed fine neighborhood $V$ of a point $p \in \mathbb C$ is called finely analytic at $p$ (on $V$) if $\mathcal{F}$ can be approximated uniformly on $V$ by analytic functions $\mathcal{F}_n$ in a neighborhood $U(\mathcal{F}_n)$ of $V$. 
\end{Def}
We will say that a continuous function on a subset $S$ of $\mathbb C$ has the Mergelyan property if it can be approximated uniformly on $S$ by analytic functions in a neighborhood of $S$. Mergelyan's Theorem states that for compact sets $K$ with finitely many components of the complement all continuous functions on $K$ which are holomorphic in the interior Int$K$ have this property.

Note that the term finely analytic functions is well known and is used for functions which have the Mergelyan property on a finely open set. Here we consider only the local definition above (we do not require that $V$ contains a fine neighborhood of each of its points). Even in this local situation a weak version of the unique continuation property holds (see Proposition 3 below). 


\begin{Def}\label{Def1}
Suppose $f$ is analytic in the unit disk $\mathbb D$. Let $p$ be a point on the unit circle $\partial \mathbb D$. We say that $f$ has fine analytic continuation $\mathcal{F}$ at $p$ if there exists a closed fine neighborhood $V$ of $p$ such that $V \cap \mathbb D \supset \overline{D(p,r)} \cap \mathbb D$ for some $r>0$, and a finely analytic function $\mathcal{F}$ at $p$ on $V$ such that $\mathcal{F}|_{\mathbb D \cap V} = f$.
\end{Def}




{\bf{Remark 1}}.
The conditions of Definition \ref{Def1} are satisfied, in particular, if $V= \overline{D(p,r)} \smallsetminus U$, where $r>0$, $U \subset \mathbb C \smallsetminus \overline{\mathbb D}$ is open and thin at $p$ and $\mathcal{F} =\mathcal{G} +\mathcal{C}$ on $V$, where $ \mathcal{G}$ is analytic on $D(p,r)$ and continuous on $\overline{D(p,r)}$, and $\mathcal{C}$ is the Cauchy-type integral of a finite Borel measure $\mu$ concentrated on $U$ such that for an increasing sequence of compacts $\kappa_n \subset U$ the functions
\begin{displaymath}
 \mathcal{C}_n(z) = -\frac{1}{\pi}\int_{\kappa_n} \frac{d\mu (\xi)}{\xi -z}, \qquad z \notin \kappa_n, 
\end{displaymath}
converge uniformly to $\mathcal{C}(z)$ on $V$.

Note that in Definition $\ref{Def1}$ we do not require that the set $V \smallsetminus \overline{\mathbb D}$ has interior points. Nevertheless, by the mentioned weak unique continuation property, if fine analytic continuation to a given set $V$ exists then it is unique on a, maybe, smaller fine neighborhood.

Here are examples of analytic functions in the unit disk which allow fine analytic continuation at certain point of the unit circle.

{\bf{Example 1}}. The functions constructed in \cite{js1} and \cite{Ediwie2} satisfy the definition. Indeed, they have the following form. Let $\{a_n\}_{n=1}^{\infty}$ be a sequence of points contained in $\mathbb C \smallsetminus \overline{\mathbb D}$ which cluster to a subset of $\partial \mathbb D$ and do not have other cluster points. Let $D(a_n, \rho_n) \subset \mathbb C \smallsetminus \overline{\mathbb D}$ be a sequence of pairwise disjoint disks around $a_n$ of radius $\rho_n >0$ such that $U\overset{def}{=} \bigcup_{n=1}^{\infty} D(a_n,\rho_n)$ is thin at a point $p \in \partial \mathbb D$. Let $c_n$ be a sequence of complex numbers such that $\sum |c_n| < \infty$ and $|c_n| \leq (1/n^2)\rho_n$. Define 
\begin{displaymath}
f(z) = \sum_{n=1}^{\infty} \frac{c_n}{z-a_n}, \qquad z \in \mathbb C \smallsetminus \bigcup_{n=1}^{\infty} D(a_n,\rho_n).
\end{displaymath}
It is immediate to check that the conditions of Definition \ref{Def1} are satisfied. \\

{\bf{Example 2}}. Let $U \subset \mathbb C \smallsetminus \overline{\mathbb D}$ be open, relatively compact and thin at every point of the unit circle $\partial \mathbb D$, and suppose moreover that $\mathbb C \smallsetminus U$ is connected. Such sets can easily be obtained by first choosing points $p_n$ which accumulates to each point of $\partial \mathbb D$ and to no other point, then choosing $\rho_n >0$ such that the disks $D(p_n,\rho_n)$ are pairwise disjoint and do not meet $\overline{\mathbb D}$. Choosing then $a_n>0$ so that the series $\sum a_n \log(|z-p_n|/ \rho_n)$ converges to a subharmonic function which is non-negative outside $\bigcup D(p_n,\rho_n)$ and, finally, choosing $r_n>0$ such that the function is less than $-1$ on $U \overset{def}{=} \bigcup D(p_n,r_n)$. Let $\mathcal{F}$ be a $C^1$ function on $\widehat{\mathbb C}$ (here $\widehat{\mathbb C}$ denotes the Riemann sphere), such that $\overline{\partial} \mathcal{F}=0$ on $\mathbb C \smallsetminus U$. By the Cauchy-Green's formula
\begin{equation*}
\mathcal{F}(z) = - \frac{1}{\pi} \iint_{U} \frac{\overline{\partial} \mathcal{F}}{\xi - z} dm_2(\xi)+ \mathcal{F}(\infty). 
\end{equation*}
Since the density $\overline{\partial} \mathcal{F}$ is bounded and the Cauchy kernel is locally integrable with respect to two-dimensional Lebesgue measure, the function $f= \mathcal{F}|_\mathbb D$ has fine analytic extension at each point $p \in \partial \mathbb D$, moreover, $\mathcal{F}$ has the Mergelyan property on $\widehat{\mathbb C} \smallsetminus U$.

More generally, let $U$ be as described and let $g \in L^{2 + \epsilon}(\mathbb C)$ for some $\epsilon>0$ and $g=0$ outside $U$. The function
\begin{displaymath}
\mathcal{F}(z) \overset{def}{=} - \frac{1}{\pi} \iint_{\mathbb C} \frac{g(\xi) dm_2(\xi)}{\xi - z},
\end{displaymath}
satisfies the condition of Definition \ref{Def1} and has the Mergelyan property on $\widehat{\mathbb C} \smallsetminus U$. If $g$ is a $C^\infty$ function then $\mathcal{F}$ is a $C^\infty$ function on the whole Riemann sphere. If $g$ is continuous and in addition $g \geq 0$ and $g>0$ at some point of each connected component of $U$, then $f = \mathcal{F}|_\mathbb D$ does not have analytic extension across any arc of $\partial \mathbb D$. Indeed, suppose it has analytic continuation $f_p$ to a disk $D(p,r)$ for some $p \in \partial \mathbb D$. By Proposition $3$ below $f_p$ coincides with the fine analytic continuation $\mathcal{F}$ on some fine neighborhood $V_1$ of $p$. $V_1$ contains a circle $\partial D(p, \rho)$, $0< \rho <r$. (This is well known. The reader who is not familiar with potential theory will find a proof below in Section 2.) Let $\kappa_n \subset U$ be an exhausting sequence of compact subsets of $U$. Since $V_1 \subset \mathbb C \smallsetminus U$ for each $n$, $\kappa_n \cap \partial D(p,\rho) = \emptyset$ and 
\begin{eqnarray*}
0 &=& \int_{|z-p|=\rho}f_p dz = \lim_{n \rightarrow \infty} - \frac{1}{\pi} \int_{|z-p|=\rho}dz \iint_{\kappa_n} \frac{\overline{\partial} \mathcal{F}}{\xi - z} dm_2(\xi) = \\
&=& \lim_{n \rightarrow \infty} - \frac{1}{\pi} \iint_{\kappa_n} dm_2(\xi) \overline{\partial}\mathcal{F} (\xi)  \int_{|z-p|=\rho} \frac{1}{\xi - z}dz = \\
&=& \lim_{n \rightarrow \infty}  \frac{2 \pi i}{\pi} \iint_{\kappa_n} dm_2(\xi) \overline{\partial}\mathcal{F} (\xi) \cdot \chi_{D(p,\rho)}(\xi) \neq 0.
\end{eqnarray*}
Here $\chi_{D(p,\rho)}$ is the characteristic function of the disk $D(p,\rho)$. The contradiction proves the assertion.

{\bf{Example 3}}. The third example is related to pseudocontinuation across certain subsets of positive length of the unit circle.
\begin{Def}
A function $f$ which is analytic in $\mathbb D$ is said to have
{\em{pseudocontinuation}} from $\mathbb D$ across the set $\mathcal{\mathcal{E}} \subset
\partial \mathbb{D}$ of positive measure to a domain $D_e \subset \{z
\in \mathbb C : |z|>1\}$ if for all $z \in \mathcal{E}$ the domain $D_e$ contains truncated non-tangential cones with vertices at $z$, and there exists an analytic
function $\widetilde{f}$ in $D_e$, such that $f$ and $\widetilde{f}$
have the same non-tangential limits at $z$. In this case we call $\widetilde{f}$ the pseudocontinuation of $f$. 
\end{Def}

For convenience we will specify the situation in the following way. We will restrict ourselves to the
case where $\mathcal{E}$ is closed and the angles and the diameters of the
truncated non-tangential cones in $D_e$ are uniformly bounded from below
by positive constants. Shrinking perhaps $\mathcal{E}$ and $D_e$ we may assume that $D_e$ is a bounded domain, moreover, that it consists of the union of all open truncated non-tangential cones with symmetry axes orthogonal
to the circle and that all the mentioned cones have the same angle and the pseudocontinuation is continuous in $\overline D_e$. So, $D_e$ has the shape of a ``saw'' near $\mathcal{E}$.  Replace ${\mathbb D}$ by a domain $D_i \subset \mathbb D$ which is symmetric to $D_e$ in a
small neighborhood of $\mathcal{E}$ in such a way that the bounded components of the complement of
$\overline{D_i} \cup \overline{D_e}$ are similar rhombs $\lozenge_l$ containing  the
complementary arcs of $\mathcal{E}$ in the circle. (The endpoints of one of the symmetry axes of the rhomb $\lozenge_l$ are the endpoints of a connected component of $\partial \mathbb D \smallsetminus \mathcal{E}$.) Assume that the unbounded component of $\mathbb C \smallsetminus ( \overline{D_e} \cup \overline{D_i})$ intersects $ \partial \mathbb D$ along a connected arc.
 We may assume that the construction is made so that the original function is continuous in
$\overline{D_i}$.

The original
function together with its pseudocontinuation across $\mathcal{E}$ defines a
continuous function on $\overline{D_i} \cup \overline{D_e}$, which
is analytic in  ${D_i} \cup {D_e}$. Denote the space of such
functions by $ A(\overline{D_i} \cup \overline{D_e})$. 

It will be convenient to give the third example with $\mathbb D$ replaced by $D_i$. It can be stated for $\mathbb D$ instead with obvious changes.

\begin{Prop}\label{Prop1}
Let $D_i,D_e$ and $\mathcal{E}$ be as above and let $f \in A(\overline{D_i} \cup \overline{D_e}) $. Put $U =\mathbb C \smallsetminus (\overline{D_i} \cup \overline{D_e})$. If $U$ is thin at a point $p \in \mathcal{E}$ and $f$ is H\"older continuous of order $\alpha \in (0,1]$, then $f|_{D_i}$ has fine analytic continuation $f|_{V_1}$ at $p$ for a fine neighborhood $V_1 \subset  \overline{D_i} \cup \overline{D_e}$ of $p$.
\end{Prop}
Note that the fact that $\mathcal{E}$ has positive length follows from the fact that its complement in $\partial \mathbb D$ is thin at $p$.

We will prove Proposition \ref{Prop1} in Section \ref{Section2}. The following Theorem holds.
\begin{Thm}\label{Thm1}
Let $f$ be analytic in $\mathbb D$ and let $p \in \partial \mathbb D$. Suppose $f$ has fine analytic continuation $\mathcal{F}$ at $p$ to a closed fine neighborhood $V$ of $p$. Then there exists another closed fine neighborhood $V_1 \subset V$ of $p$, such that the graph $\Gamma_{\mathcal{F}} (V_1)$ is contained in the pluripolar hull of $\Gamma_{f}(\mathbb D)$. 

Moreover, if Int$V \smallsetminus \overline{\mathbb D}$ has a connected component $\overset\circ V$ which is not thin at $p$ then $\Gamma_{\mathcal{F}} (\overset\circ V)$ is contained in the pluripolar hull of $\Gamma_f(\mathbb D)$. 

If $V= \overline{D(p,r)}\smallsetminus U$ with $U$ open and thin at $p$ and $\overline{U}\smallsetminus \overline{\mathbb D}$ is also thin at $p$ then there exists a unique connected component of Int$V \smallsetminus \overline{\mathbb D} = \{ |z-p| < r, |z| > 1 \} \smallsetminus \overline{U}$ which is not thin at $p$.
\end{Thm}
Note that Theorem \ref{Thm1} holds as well in the situation of Proposition \ref{Prop1} with $\mathbb D$ replaced by $D_i$.
Theorem \ref{Thm1} can be slightly generalized.
\begin{Thm}\label{Thm1a}
Let $\mathcal{F}$ be finely analytic on a closed fine neighborhood $V= \overline{D(p,r)} \smallsetminus U$ of a point $p \in \mathbb C$. Let $\gamma$ be a smooth arc, $\gamma: [-1,1] \rightarrow \mathbb C$ with $\gamma(0) =p$, which divides $D(p,r)$ into two components $D_+(p,r)$ and $D_-(p,r)$. Suppose $\overline{U} \smallsetminus \gamma$ is thin at $p$. Then there are unique connected components $V_+$ and $V_-$ of $D_+(p,r) \smallsetminus \overline{U}$ and $D_-(p,r) \smallsetminus \overline{U}$ which are not thin at $p$. For those we have that the pluripolar hull of $\Gamma_{\mathcal{F}}(V_+)$ is contained in $\Gamma_{\mathcal{F}}(V_-)$ and vice versa.
\end{Thm}

\begin{Cor}\label{Cor1}
Let $D_i,D_e$ and $\mathcal{E}$ be as described before Proposition \ref{Prop1} and let $f\in A(\overline{D_i} \cup \overline{D_e})$. If $U = \mathbb C \smallsetminus (\overline{D_i} \cup \overline{D_e})$ is thin at some point $p \in \mathcal{E}$ and $f$ is H\"older continuous of order $\alpha \in (0,1]$ then $\Gamma_f(D_e)$ is contained in the pluripolar hull of $\Gamma_f(D_i)$. 
\end{Cor}
Theorem \ref{Thm1a} and Corollary \ref{Cor1} has the following consequence.
\begin{Cor}\label{Cor2}
There exist univalent analytic functions in the disk which are smooth up to the boundary and are nowhere analytically continuable and have an analytic manifold in the
non-trivial part of the pluripolar hull of the graph.
\end{Cor}

Functions with the mentioned property except univalency were constructed in \cite{Ediwie2}. Theorem \ref{Thm1} and Example $2$ give further classes of functions with the mentioned property (not necessarily univalent functions). The present constructions are simpler then those in \cite{Ediwie2}.

Denote by $\pi_1$ the projection onto the first coordinate plane in $\mathbb C^2$, $\pi_1(z_1,z_2)= z_1$ for $z =(z_1,z_2) \in \mathbb C^2$. Theorem \ref{Thm1} and its corollaries have the following counterpart.
\begin{Thm}\label{Thm2}
Let $f$ be analytic in $\mathbb D$. Suppose $(\Gamma_f(\mathbb D))^*_{\mathbb C^2}$ is not contained in $\mathbb D \times \mathbb C$. Put $E\overset{def}{=} \overline{(\Gamma_f(\mathbb D))^*_{\mathbb C^2}}.$ Then $\pi_1(E)$ contains a fine neighborhood of a point $p \in \partial \mathbb D$ (i.e. $\mathbb C \smallsetminus \pi_1(E)$ is thin at p).
\end{Thm}
For bounded analytic functions the following sharper version of Theorem \ref{Thm2} holds.
\begin{Thm}\label{Thm3a}
Let $f$ be a bounded analytic function in $\mathbb D$ whose graph is not complete pluripolar in $\mathbb C^2$. Then $\pi_1(\overline{(\Gamma_f(\mathbb D))^*_{\mathbb C^2}})$ contains a fine neighborhood of a point $p \in \partial \mathbb D$.
\end{Thm}
The reduction of Theorem \ref{Thm3a} to Theorem \ref{Thm2} follows immediately from results in \cite{Ediwie3}.

Lemma \ref{Lem2} in Section $3$ shows that Theorem \ref{Thm2} can be slightly improved. We will not give the corresponding statements here.
\begin{Cor}\label{Cor3}
Let $\mathcal{D} = \mathbb D \cup D$, where $D$ is a domain, $D \subset \mathbb C \smallsetminus  \overline{\mathbb D}$, which contains a truncated non-tangential cone of fixed size with vertex $z$ for each $z \in \mathcal{E}$, $\mathcal{E}$ a closed subset of positive length on $\partial \mathbb D$. Let $f$ be holomorphic in $\mathcal{D}$ and continuous in $\overline{\mathcal{D}} = \overline{\mathbb D}\cup \overline{D}$ (hence $f|_{\mathbb D}$ and $f|_{D}$ are pseudocontinuations of each other across the set $\mathcal{E}$). Suppose $\mathbb C \smallsetminus \overline{\mathcal{D}}$ is non-thin at every point $z\in \partial \mathbb D$ and $\Gamma_f(\overline{\mathcal{D}})$ is complete pluripolar. Then $\Gamma_f(\mathbb D)$ is complete pluripolar.
\end{Cor}
Note that functions $f$ with the mentioned properties exist (see \cite{Edl}, to appear). Corollary \ref{Cor3} states roughly that if two analytic manifolds in $\mathbb C^2$ have contact along a set which is not massive enough in potential theoretic terms then the property of plurisubharmonic functions in $\mathbb C^2$ to be $- \infty$ does not propagate from one of the manifolds to the other one.
\begin{proof}[Proof of Corollary \ref{Cor3}]
$E= \overline{(\Gamma_f(\mathbb D)^*_{\mathbb C^2})} \subset \Gamma_f(\overline{\mathcal{D}})$ since the latter set is closed and complete pluripolar. Hence $\pi_1(E) \subset \overline{\mathcal{D}}$ and $\mathbb C \smallsetminus \pi_1(E) \supset \mathbb C \smallsetminus \overline{\mathcal{D}}$ is non-thin at any point $p \in \partial \mathbb D$.
\end{proof} 

Note that the set $V$ in Theorem \ref{Thm1} and the set $E$ in Theorem \ref{Thm2} were assumed to be closed. It would be interesting to remove the condition of closeness.
We conclude this Section with an example of fine analytic continuation to a set with no interior points outside the unit disk, with an example with infinitely sheeted pluripolar hull and with some open problems.

{\bf{Problem 1}}. Is $\pi_1((\Gamma_f(\mathbb D))^*_{\mathbb C^2})$ finely open ?

{\bf{Problem 2}}. Let $f$ be analytic in $\mathbb D$, not necessarily bounded. Suppose that the pluripolar hull of $\Gamma_f(\mathbb D)$ is contained in $\mathbb D \times \mathbb C$. Does this imply that $\Gamma_f(\mathbb D)$ is complete pluripolar in $\mathbb C^2$ ?

{\bf{Example 4}}. Consider all points in $\mathbb C \smallsetminus
\overline{\mathbb D}$ with rational coordinates. This set is countable
and accumulates, in particular to the whole circle $\partial \mathbb
D$. As in Example 2 there exists a subharmonic function $\mathcal{U}$
which is $-\infty$ on this set and non-negative at each point $p \in
\overline{\mathbb D}$. Let $U$ be the set of points on which
$\mathcal{U} < -1$. $U$ is open, contained in $\mathbb C
\smallsetminus \overline{\mathbb D}$ and it is thin at each point of
the unit circle. Moreover $\mathbb C \smallsetminus U$ is connected,
i.e. $U$ is simply connected which is a consequence of the maximum principle.

Let $\mathcal{F}$ be a continuous function on $\mathbb C \smallsetminus U$ which has the Mergelyan property on each compact subset of $\mathbb C \smallsetminus U$ and has complete pluripolar graph $\Gamma_{\mathcal{F}}(\mathbb C \smallsetminus U)$. Such functions were constructed in \cite{Edl} for arbitrary closed subsets of $\mathbb C$. Then $f = \mathcal{F}|_\mathbb D$ is analytic and $f$ has fine analytic continuation at each point $p \in \partial \mathbb D$. Hence by Theorem \ref{Thm1a} there exists a set $A_1$ which contains a fine neighborhood of each point of the circle, such that the graph $\Gamma_\mathcal{F}(A_1)$ is in the pluripolar hull of $\Gamma_\mathcal{F}(\mathbb D)$. However the mentioned pluripolar hull is contained in $\Gamma_\mathcal{F}(\mathbb C \smallsetminus U)$, and the set $\mathbb C \smallsetminus U$ does not have interior points outside the unit disk $\mathbb D$. Hence in this case the non-trivial part of the pluripolar hull of $\Gamma_{f}(\mathbb D)$ does not have analytic structure, i.e. there is no piece of a non-trivial analytic manifold contained in it. The following problems arise.

{\bf{Problem 3}}. Let $f$ be analytic in $\mathbb D$. Suppose $\pi_1((\Gamma_f(\mathbb D))^*_{\mathbb C^2}$ has an interior point $p$ outside $\mathbb D$. (Note that here we consider the pluripolar hull itself, not its closure!) Is there a neighborhood $V_p$ of $p$ in $\mathbb C$ and a relatively closed subset $X$ of $V_p \times \mathbb C$ such that $X \subset (\Gamma_f(\mathbb D))^*_{\mathbb C^2}$ and $X$ is a ``limit'' (e.g. in the Hausdorff topology) of relatively closed analytic varieties in $V_p \times \mathbb C$?

{\bf{Problem 4}}. Let $f$ be analytic in $\mathbb D$ and suppose that $\Gamma_f(\mathbb D)$ is not complete pluripolar. How big can the fiber of  $(\Gamma_f(\mathbb D))^*_{\mathbb C^2}$ be over a ``generic point'' in $\pi_1((\Gamma_f(\mathbb D))^*_{\mathbb C^2})$? Can it be more than countable ? Does $(\Gamma_f(\mathbb D))^*_{\mathbb C^2}$ contain the graph of a reasonable function over a sufficiently massive subset of $\mathbb C$ which reflects certain generalized analytic continuation property of $f$ ? 

We have the following example in mind.

{\bf{Example 5}}. By a segment we mean the (closed) part of a real line in $\mathbb C$ which joins two points in $\mathbb C$. Let $U$ be as in Example $2$. Consider a sequence of pairwise disjoint segments $\sigma_n$ contained in $U$ which accumulates to the whole circle $\partial \mathbb D$ and does not have other limit points in $\mathbb C$. 
Denote the endpoints of the segment $\sigma_n$ by $a_n$ and $b_n$ and associate to $\sigma_n$ the branch of the function $\sqrt{(z-a_n)(z-b_n)}$ on $\mathbb C \smallsetminus \sigma_n$ which equals $z+O(1)$ near $\infty$. We will use this notation only for the mentioned branch. Let $c_n$ be complex numbers so that $\sum |c_n| < \infty$. Let $\mathcal{F}(z)= \sum c_n \sqrt{(z-a_n)(z-b_n)}$ on $\mathbb C \smallsetminus U$. The series converge uniformly on compacts in $\mathbb C \smallsetminus U$. The function $f = \mathcal{F}|_{\mathbb D}$ has fine analytic continuation at each point of $\partial \mathbb D$. Moreover, $\mathcal{F}$ has the Mergelyan property on compacts in $\mathbb C \smallsetminus U$. By Theorem \ref{Thm1} the graph $\Gamma_{\mathcal{F}}(\{|z|>1\} \smallsetminus \overline{U})$ is in the pluripolar hull of $\Gamma_f(\mathbb D)$. Since $\mathcal{F}$ has single-valued analytic continuation to $\{|z|>1\} \smallsetminus \bigcup \sigma_n$, the graph of $\mathcal{F}$ over this set is contained in the pluripolar hull of $\Gamma_f(\mathbb D)$ too.

Fix a number $n$. The graph of the function $\sqrt{(z-a_n)(z-b_n)}$ over $\mathbb C \smallsetminus \sigma_n$ is an open subset of the algebraic curve $\{(z,w) \in \mathbb C^2: w^2 = (z-a_n)(z-b_n) \}$. There is a neighborhood $U_n$ of $\sigma_n$ such that $\sum_{l \neq n}c_l \sqrt{(z-a_l)(z-b_l)}$ is analytic in $U_n$. Hence the analytic set
\begin{displaymath}
A_n = \{ (z,w) \in U_n \times \mathbb C : \big ( w- \sum_{l \neq n}c_l \sqrt{(z-a_l)(z-b_l)} \big )^2 = c_n^2 (z-a_n)(z-b_n) \}
\end{displaymath}
contains the graph $\Gamma_{\mathcal{F}}(U_n \smallsetminus \sigma_n)$ and is therefore in the pluripolar hull of $\Gamma_{\mathcal{F}}(U_n \smallsetminus \sigma_n)$ and, hence, of $\Gamma_f(\mathbb D)$. The set $A_n$ has two sheets over $U_n \smallsetminus \sigma_n$, the second sheet is the graph of the function 
\begin{displaymath}
-c_n\sqrt{(z-a_n)(z-b_n)} + \sum_{l \neq n}c_l \sqrt{(z-a_l)(z-b_l)}
\end{displaymath} 
This function has analytic continuation to $(\mathbb C \smallsetminus \overline{\mathbb D}) \smallsetminus \overline{U}$. Moreover, it has the Mergelyan property on compact subsets of $\mathbb C \smallsetminus U$. By Theorem \ref{Thm1a} it has fine analytic continuation at each point of $\partial \mathbb D$ to the disk $\mathbb D$. We obtained that the pluripolar hull of $\Gamma_f(\mathbb D)$ contains a two-sheeted branched covering over the set $\mathbb D \cup (\mathbb C \smallsetminus (\overline{\mathbb D} \cup\bigcup_{l \neq n} \sigma_l))$, the points $a_n$ and $b_n$ being the branch points. Repeating the argument for all other endpoints of the segments, we obtain that the pluripolar hull of $\Gamma_f(\mathbb D)$ contains an infinitely-sheeted branched covering (countably many sheets over generic points) over the set $\mathbb C \smallsetminus \partial \mathbb D$ with branch points $\{a_n\}$ and $\{b_n\}$. Note that the sheets over $\mathbb D$ are graphs of analytic functions. Moreover, the covering is unbranched over $\mathbb C \smallsetminus ( \mathbb T \cup \bigcup_n(\{a_n\}\cup \{b_n\}))$. The infinitely-sheeted branched covering over $\mathbb C \smallsetminus \partial \mathbb D$ can be approximated by analytic subsets of $(\mathbb C \smallsetminus \partial \mathbb D) \times \mathbb C$, the sheets of which over $\mathbb C \smallsetminus \bigcup_1^n \sigma_l$ are the graphs of $\mathcal{F}_n(z) = \sum_{l=1}^n \pm c_l \sqrt{(z-a_l)(z-b_l)}$ with all possible choices of $+$ and $-$. Note that similar arguments as in Example $2$ show that the function $f$ is nowhere analytically extendible across $\partial \mathbb D$. Choosing the $c_n$ more carefully we may reach that $f$ is of class $C^\infty$ on the closed unit disk.

When this paper was written we received a preprint of Zwonek \cite{zw} where he constructed an analytic function in $\mathbb D$ which does not have analytic extension across $\partial \mathbb D$ and for which $(\Gamma_f(\mathbb D))^*_{\mathbb C^2}$ has at least two sheets over $\mathbb D$. 

{\bf{Problem 5}}. Is $(\Gamma_f(\mathbb D))^*_{\mathbb C^2}$ related to a suitable positive $(1,1)$-current ?

In Section \ref{Section2} of the paper we will prove Proposition \ref{Prop1}, Theorem \ref{Thm1} and its corollaries. In the remaining Section \ref{Section3} we will prove Theorem \ref{Thm2}. 

\section{Non-Trivial hull.}\label{Section2}

In this section we will prove Proposition \ref{Prop1}, Theorem $\ref{Thm1}$ and its corollaries. Recall that in Proposition \ref{Prop1} we consider domains $D_i$ and $D_e$, $D_i \subset \mathbb D$ and $D_e \subset \mathbb C \smallsetminus \overline{\mathbb D}$, such that the bounded components of $\mathbb C \smallsetminus (\overline{D_i} \cup \overline{D_e})$ are similar rhombs $\lozenge_l$ for which the endpoints of one of the symmetry axis are the endpoints of a connected component of $\partial \mathbb D \smallsetminus \mathcal{E}$. Here $\mathcal{E} \overset{def}{=}  \overline{D_i} \cap \overline{D_e}$. The following Lemma will be useful.
\begin{Lem}\label{Lem1}
Let $\bigcup_l I_l$ be a union of open pairwise disjoint arcs on $\partial \mathbb D$ which is thin at $p \in \partial \mathbb D$. Denote by $\bigcup_l \overline{\lozenge}_l$ the union of closed similar rhombs with the property that two opposite corners of $\overline{\lozenge}_l$ are the endpoints of $I_l$. Then $\bigcup \overline{\lozenge}_l$ is thin at $p$. 
\end{Lem}
For the proof we need the following proposition which is interesting in itself.
\begin{Prop}\label{thinp}
Let $E \subset \mathbb C$ be thin at the point $0 \in \overline{E}$. Let $T: \mathbb C \rightarrow \mathbb C$ be a mapping which satisfies a Lipschitz condition ($|Tz_1 - Tz_2| \leq C|z_1-z_2|$ for $z_1,z_2 \in \mathbb C$) and such that $T(0)=0$ and $|Tz| \geq c|z|$ for $z \in \mathbb C$. ($C$ and $c$ are positive constants). Then the set $TE$ is thin at 0.
\end{Prop}
This proposition was known already to Brelot. Since a slightly weaker assertion is stated in \cite{Bre} (see chapter 7, paragraph 2 ) and we were not able to find an explicit reference for Proposition \ref{thinp}, we sketch the proof for the convenience of the reader. 

\begin{proof}
Since $E$ is thin at $0$ there exists a subharmonic function $\mathcal{V}$ in a neighborhood of $0$ with $\mathcal{V}(0) > - \infty$ and $\lim_{\xi \in E, \xi \rightarrow 0} \mathcal{V}(\xi) = - \infty$. Using the Riesz representation theorem and subtracting a harmonic function we may assume that $\mathcal{V}$ has the form 
\begin{displaymath}
\mathcal{V}(z) = \int \log| \xi - z| d \mu (\xi)
\end{displaymath}
for a positive Borel measure $\mu$. Define the measure $\mu_1$, $\mu_1(A) = \mu (T^{-1}(A))$ for each Borel set $A$, and put 
\begin{displaymath}
\mathcal{V}_1(z) = \int \log| \xi - z| d \mu_1 (\xi).
\end{displaymath}
Then 
\begin{displaymath}
\mathcal{V}_1(Tz) =  \int \log| \xi - Tz| d \mu_1 (\xi) =  \int \log| T\xi - Tz| d \mu (\xi).
\end{displaymath}
Hence $\mathcal{V}_1(Tz) \leq \mathcal{V}(z) + \log C \cdot || \mu ||$ and $\mathcal{V}_1(0) \geq \mathcal{V}(0) + \log c \cdot || \mu ||$.
\end{proof}

\begin{proof}[Proof of Lemma \ref{Lem1}]
Note first that $\bigcup_l \overline{I}_l$ is thin at $p$ if  $\bigcup_l I_l$ is, since the sets differ by a countable (and hence thin) set. Denote by $\Lambda$ the union of the boundaries of the rhombs. Let $\Lambda_+$ and $\Lambda_-$ be the parts of $\Lambda$ which are contained outside the unit disk and inside the closed unit disk, respectively. Both $\Lambda_+$ and $\Lambda_-$, can be represented as graphs over a part of $\partial \mathbb D$; 
\begin{displaymath}
\Lambda_+ =\{ r_+e^{i \phi} | \textrm{  } r_+= r_+(\phi),\textrm{  } e^{ i \phi} \in \bigcup_l \overline{I_l} \} \textrm{ and } \Lambda_- =\{ r_-e^{i \phi} |\textrm{  } r_-= r_-(\phi),\textrm{  } e^{ i \phi} \in \bigcup_l \overline{I_l} \}.
\end{displaymath}
The mapping $T(e^{i \phi })=r_+ e^{i \phi}$ where $e^{i \phi} \in \bigcup_l \overline{I_l}$ can be extended to the whole plane as a Lipschitz continuous mapping which satisfies the conditions of Proposition \ref{thinp} with $0$ replaced by $p$. (The same is true for the mapping $T(e^{i \phi })=r_- e^{i \phi}$). Since thinness is invariant under translation, Proposition \ref{thinp} shows that both, $\Lambda_+$ and $\Lambda_-$, are thin at $p$. Therefore their union $\Lambda_+ \cup \Lambda_- = \Lambda$ is thin at $p$. Since the union of the boundaries of the closed rhombs is thin at $p$ we conclude that the union of the closed rhombs is thin at $p$.
\end{proof}
We need the following immediate consequence of Lemma \ref{Lem1}.
\begin{Cor}\label{Cor4}
If $U = \bigcup_l \lozenge_l$ is thin at $p \in \partial \mathbb D$ then also $\overline{U} \smallsetminus \overline{\mathbb D}$ is thin at $p$.
\end{Cor}
\begin{proof}
Indeed, $\overline{U} \smallsetminus \overline{\mathbb D} = \bigcup_l \overline{ \lozenge_l}\smallsetminus \overline{\mathbb D}$.
\end{proof}
We will make use also of the following three observations. 

If the union of rhombs $\bigcup \overline{ \lozenge_l}$ is thin at $p$ there are arbitrarily small numbers $r >0$ with the property that $\partial D(p,r) \cap \bigcup \overline{ \lozenge_l} = \emptyset$. (See \cite{Bre} or Proposition \ref{thinp} with $Tz =|z|$ after suitable translation).

Moreover, looking at connected components of the union of closed intervals and replacing the previous rhombs by
closed rhombs associated to these connected components in the same way as
above, we may assume that the $\overline{\lozenge}_l$:s are pairwise disjoint.
 
If $\bigcup \overline{ \lozenge_l}$ is thin at $p$ then there exists another sequence of similar (open) rhombs $\lozenge_j'$ associated to disjoint open arcs $I_j'$ of $\partial \mathbb D$, such that   $\bigcup \lozenge_j'$ is thin at $p$ and $\bigcup \overline{ \lozenge_l} \subset \bigcup \lozenge_j'$.
In fact, $\bigcup \overline{I}_l$ is thin at $p$. If for a subharmonic function $\mathcal{V}$ $\overline{\lim}_{z \in \bigcup \overline{I}_l, z \rightarrow p} \mathcal{V}(z) < a < \mathcal{V}(p)$, then for each $\overline{I}_l$ contained in a small neighborhood $V_p$ of $p$ $\sup_{\overline{I}_l} \mathcal{V} < a$, hence $\sup_{\widetilde{I}_l} \mathcal{V} < a$ for some open arc $\widetilde{I}_l \supset \overline{I}_l$. Take also for the other arcs $\overline{I}_l$ suitable open arcs $\widetilde{I}_l$ containing them. The set $ \bigcup_l \widetilde{I}_l$ is thin at $p$. Let the $I_j'$ be the connected components of the latter union and associate rhombs $\lozenge_j'$ to them.
\begin{proof}[Proof of Proposition \ref{Prop1}]
We will assume that the $\lozenge_l$ in the statement of Proposition \ref{Prop1} are pairwise disjoint (shrinking otherwise the sets $D_i$ and $D_e$) and prove fine analytic continuation to $\overline{D(p,r)} \smallsetminus \bigcup_j \lozenge_j'$ for a suitable small $r>0$ and the rhombs $\lozenge_j'$ described above. We may assume that $r>0$ is chosen so that $\partial D(p,r)$ does not meet $\bigcup_j \lozenge_j'$ (by the remark above). Keep notation $\overline{ \lozenge_l}$ for only those of the original rhombs which are contained in $D(p,r)$ and $\lozenge_j'$ for those of the enlarged rhombs which are contained there. The function $f$ is analytic in each of the domains $\mathcal{D}_i' \overset{def}{=} D(p,r) \cap \mathbb D \smallsetminus \bigcup_l \overline{\lozenge_l}$ and $\mathcal{D}_e' \overset{def}{=} D(p,r) \cap (\mathbb C \smallsetminus \overline{\mathbb D}) \smallsetminus \bigcup_l \overline{\lozenge_l}$ and H\"older continuous in the union of the closures. Both domains have rectifiable boundary, hence by Cauchy's formula
\begin{displaymath}
f(z)= \frac{1}{2 \pi i}\int_{\partial \mathcal{D}_i'} \frac{f(\xi)}{\xi - z}d \xi + \frac{1}{2 \pi i}\int_{\partial \mathcal{D}_e'} \frac{f(\xi)}{\xi - z}d \xi, \qquad z \in \mathcal{D}_i' \cup \mathcal{D}_e'.
\end{displaymath}
The contours of integration are always oriented as boundaries of relatively compact domains. Note that one of the integrals in the sum above will be equal to zero. Using that $\partial \mathbb D \cap  D(p,r) \smallsetminus \bigcup \lozenge_l = \mathcal{E} \cap D(p,r)=\partial \mathcal{D}_i'\cap \partial\mathcal{D}_e'$ and orientation over this set is provided twice with opposite orientation we obtain
\begin{eqnarray*}
f(z) &=& \frac{1}{2 \pi i}\int_{\partial D(p,r)} \frac{f(\xi)}{\xi - z}d \xi - \frac{1}{2 \pi i} \sum_l \int_{\partial \lozenge_l} \frac{f(\xi)}{\xi - z}d \xi \overset{def}{=} \\
&\overset{def}{=}& J(z) - \sum_l J_l(z), \qquad z \in \mathcal{D}_i' \cup \mathcal{D}_e'.
\end{eqnarray*}
By Privalov's theorem $J(z)$ extends to a H\"older continuous function of order $\alpha$ in $\overline{D(p,r)}$ if $\alpha <1$ and of any order less than $1$ if $\alpha =1$. The measure $f(\xi) d \xi$ on $\bigcup \partial \lozenge_l$ is a finite Borel measure concentrated on a subset of $\bigcup \lozenge_j'$. To prove Proposition \ref{Prop1} let $\kappa_n = \bigcup_{l=1}^n \overline{\lozenge_l}$ and $\mathcal{F}_n(z)= J((1-1/n)z) - \sum_{l=1}^n J_l(z)$, $ z \in  \overline{D(p,r)} \smallsetminus \kappa_n$. We have to check that the  $\mathcal{F}_n$ converge uniformly to $f$ on $\overline{\mathcal{D}}_i' \cup \overline{\mathcal{D}}_e' = \overline{D(p,r)} \smallsetminus \bigcup_j \lozenge_j'$. To obtain a uniform estimate of $J_l$ on $ \overline{D(p,r)} \smallsetminus \lozenge_l$ we use that for $z \notin \overline{\lozenge_l}$ the Cauchy type integral of the constant function $f(z)$ with pole at $z$ along $\partial \lozenge_l$ vanishes. We get for $z \in  \overline{D(p,r)} \smallsetminus \overline{\lozenge_l}$
\begin{displaymath}
|J_l(z)|= |\frac{1}{2 \pi i}\int_{\partial \lozenge_l } \frac{f(\xi)- f(z)}{\xi - z}d \xi| \leq 
 C  \int_{\partial \lozenge_l} \frac{|\xi -z |^{\alpha}}{|\xi -z |}|d \xi |.
\end{displaymath}
Let $N>n$ be a natural number and let $z \in  \overline{D(p,r)} \smallsetminus \bigcup_{l=n}^N \overline{\lozenge_l}$. Then
\begin{displaymath}
\sum_{n \leq l \leq N} |J_l(z)| \leq C \int_{\bigcup_{n \leq l \leq N} \partial \lozenge_l} \frac{|\xi -z |^{\alpha}}{|\xi -z |}|d \xi |   \leq C \int_{\bigcup_{l \geq n} \partial \lozenge_l} \frac{|\xi -z |^{\alpha}}{|\xi -z |}|d \xi |.
\end{displaymath}
Represent the contour of integration on the right hand side as the union of its part $\Lambda_k^-$ contained in $\overline{\mathbb D}$ and its part $\Lambda_k^+$ contained in $\mathbb C \smallsetminus \overline{\mathbb D}$. Each of the parts is the graph of a Lipschitz continuous function over a subset of $(\partial \mathbb D) \cap \overline{D(p,r)}$ with uniform estimate for the Lipschitz constant (which depends on the angle of the truncated non-tangential cones contributing to $D_i$ and $D_e$). For a point $\zeta \neq 0$ we denote by $\zeta'$ its radial projection to the circle, $\zeta' = \zeta / |\zeta|.$ Using the inequality $| \xi - z | \geq \textrm{const} | \xi' -z'|$ and estimating the arc-length on $\Lambda_k^+$ and on $\Lambda_k^-$ by arc-length of the radial projection we obtain

\begin{eqnarray*}
\sum_{n \leq l \leq N} |J_l(z)| &\leq& C \int_{\partial \mathbb D  \cap \bigcup_{l \geq n} \overline{\lozenge_l}}|\xi' -z' |^{\alpha - 1}|d \xi' | \leq  \\
&\leq& C' \int_{\gamma_n}|e^{i \phi} -1|^{\alpha - 1}|d e^{i \phi}|,\qquad z \in \overline{D(p,r)}\smallsetminus \bigcup_{l=n}^N \overline{\lozenge_l}.
\end{eqnarray*}
where $\gamma_n$ is the arc of the circle which is symmetric around the point $1$ and has length mes$_1(\partial \mathbb D \cap \bigcup_{l \geq n}\overline{\lozenge_l})$. Since $\alpha > 0$ the right hand side converge to zero for $n \rightarrow \infty$. This proves the proposition.
\end{proof} 

For a domain $G \subset \mathbb C$, a Borel subset $\mathcal{E}$ of $\partial G$ and a point $z \in G$ we denote by $\omega(z,\mathcal{E},G)$ the harmonic measure of $\mathcal{E}$ with respect to $G$ computed at the point $z$.
\begin{proof}[Proof of Theorem \ref{Thm1}]
Let $V = \overline{D(p,r)} \smallsetminus U$, where $U \subset \mathbb C \smallsetminus \overline{\mathbb D}$ is open and thin at $p$. We will first obtain a harmonic measure estimate. Let $\rho > 0$ be small enough and such that $\{|z-p| = \rho \} \cap U = \emptyset$. Since $U$ is thin at $p$ such $\rho$ exists. Let $J$ be a closed subarc of $\partial D(p,\rho)$ contained in $\mathbb D$. Decreasing $\rho$ we may assume that the length of $J$ is at least $5\pi \rho /6$. Let $K_n$ be an increasing sequence of compact subsets of $U$, each $K_n$ being the finite disjoint union of closures of simply connected domains. Then $D(p,r) \smallsetminus K_n$ is connected.
We claim that if $\rho$ is small enough there exists a number $r_1$, $ 0< r_1 < \rho$ and an open set $U_1 \supset U \cap D(p, \rho)$, which is thin at $p$ such that the following harmonic measure estimate holds:
\begin{equation}\label{hm1}
\omega(z,J, D(p,\rho) \smallsetminus K_n) \geq \frac{1}{4} \textrm{ for each } z \in V_1= D(p,r_1) \smallsetminus U_1 \textrm{ and each } n.
\end{equation} 
In fact, since $U$ is thin at $p$ there is a subharmonic function $\mathcal{U}$ in a neighborhood of $p$ which is finite at $p$ and tends to $-\infty$ along the set $U$. Taking $\rho$ small enough and adding a constant to $\mathcal{U}$ we may assume that $\mathcal{U}$ is defined and $<0$ in $\overline{D(p,\rho)}$. Multiplying $\mathcal{U}$ by a positive constant we may assume that $\mathcal{U}(p) > -1/12$. Taking $\rho$ small enough, we may assume that $\mathcal{U} < -1$ on $U \cap \overline{D(p,\rho)}$. Then
\begin{equation}\label{hm2}
\omega(z,J, D(p,\rho) \smallsetminus K_n) \geq \omega(z,J, D(p,\rho)) + \mathcal{U}(z), \qquad z \in D(p,\rho) \smallsetminus K_n.
\end{equation} 
Indeed, the boundary of $D(p,\rho) \smallsetminus K_n$ is smooth, hence regular for the Dirichlet problem. The left hand side is harmonic in this domain and extends continuously to all but two points of its closure, the right hand side is subharmonic in the domain, bounded from above and its boundary values at all but two points are majorized by those on the left hand side. Denote by $U_1$ the set $U_1\overset{def}{=} \{ z \in D(p, \rho): \mathcal{U}(z) < -1/12 \}$. $U_1$ is open and since $\mathcal{U}(p) > -1/12$ the set $U_1$ is thin at $p$. Clearly  $U_1 \supset U\cap D(p, \rho )$. By the assumption on the length of $J$ we have $\omega(p,J, D(p,\rho)) \geq 5/12$. Let $r_1 \in (0,\rho)$ be so small so that 
\begin{equation}\label{hm3}
\omega(z,J, D(p,\rho)) > \frac{4}{12} \textrm{ for } z \in D(p,r_1).
\end{equation} 
Then by (\ref{hm2}), the definition of $U_1$ and by (\ref{hm3})
\begin{equation*}\label{hm4}
\omega(z,J, D(p,\rho) \smallsetminus K_n) > \frac{4}{12}- \frac{1}{12} = \frac{1}{4} \textrm{ for } z \in D(p,r_1) \smallsetminus U_1 \overset{def}{=}V_1 \textrm{ for each } n.
\end{equation*} 
(\ref{hm1}) is proved.

Suppose now that $f$ has fine analytic continuation $\mathcal{F}$ at $p$ to a fine neighborhood $V = \overline{D(p,r)} \smallsetminus U$, i.e. there exist analytic functions $\mathcal{F}_n$ in neighborhoods $U(\mathcal{F}_n)$ of $V$ which converge uniformly to $\mathcal{F}$ on $V$. Shrinking the neighborhoods $U(\mathcal{F}_n)$ we may always assume that $\sup_{U(\mathcal{F}_n)}|\mathcal{F}_n| \leq C$ for all $n$ and a constant $C>1$. Since $U$ is simply connected one can choose an increasing sequence of compact subsets $K_n$ of $U$, each being the finite disjoint union of closures of simply connected domains such that $\overline{D(p,r) \smallsetminus K_n} \subset U(\mathcal{F}_n)$. Hence $\mathcal{F}_n$ are analytic in $D(p,r) \smallsetminus K_n$ and continuous in $\overline{D(p,r) \smallsetminus K_n}$ and their maximum norms in these sets are bounded by the constant $C$. Take $\rho$, $r_1$ and $V_1$ as above. Fix an arbitrary point $z \in V_1$ and define $Q_{n,z}(\xi) = \mathcal{F}_n(\xi) +\mathcal{F}(z) - \mathcal{F}_n(z)$, $\xi \in \overline{D(p,r) \smallsetminus K_n}$. Then $Q_{n,z}$ are analytic and uniformly bounded on $D(p,\rho) \smallsetminus K_n$ and continuous on $\overline{D(p,\rho) \smallsetminus K_n}$, $Q_{n,z}(z) = \mathcal{F}(z)$, and $Q_{n,z}\rightarrow \mathcal{F}$ uniformly on $V$ for $n \rightarrow \infty$. 

Let $B$ be a ball in $\mathbb C^2$ which contains the graphs $\Gamma_{Q_{n,z}}(D(p,r) \smallsetminus K_n)$ for all $n$ and $z$. Fix $z \in V_1$. Let $u$ be a plurisubharmonic function in $\mathbb C^2$ which equals $- \infty$ on $\Gamma_f(\mathbb D)$. Adding a constant, we may assume that $u<0$ in $B$. $J \subset V \cap \mathbb D$ and for large $n$ the set $\Gamma_{Q_{n,z}}(J)$ is uniformly close to $\Gamma_f(J) = \Gamma_\mathcal{F}(J)$. Since $u$ is upper semi-continuous for each large $N$ there exists $n$ such that $u(\xi,Q_{n,z}(\xi)) < -N$ for $\xi \in J$. The function $\xi \mapsto  u(\xi,Q_{n,z}(\xi))$ is a negative subharmonic function on $D(p, \rho) \smallsetminus K_n$ which is  upper semi-continuous on the closure of this set, hence by the estimate of harmonic measure we obtain
\begin{equation}\label{hm5}
u(z,\mathcal{F}(z)) = u(z, Q_{n,z}(z)) \leq -N \omega(z,J, D(p,\rho) \smallsetminus K_n) < -\frac{N}{4}.
\end{equation} 
Since $N$ was arbitrary we obtain $u(z,\mathcal{F}(z)) = - \infty$ for all $z \in V_1$ if $u = - \infty$ on $\Gamma_f(\mathbb D)$.

Suppose now that Int$V \smallsetminus \overline{\mathbb D}$ has a connected component $\overset\circ V$ which is non-thin at $p$. Then $\overset\circ V \cap V_1 = \overset\circ V \cap \overline{D(p,r_1)} \smallsetminus U_1$ is non-thin at $p$ (since $\overset\circ V \subset ( \overset\circ V \cap  \overline{D(p,r_1)} \smallsetminus U_1) \cup \{|z-p| > r_1 \}  \cup U_1$ and the last two sets are thin at $p$). Hence $\overset\circ V \cap V_1$ is not polar, and therefore, since $\mathcal{F}$ is analytic on $\overset\circ V$, $\Gamma_\mathcal{F}(\overset\circ V)$ is contained in the pluripolar hull of $\Gamma_f(\mathbb D)$.

Suppose $\overline{U} \smallsetminus \overline{\mathbb D}$ is thin at $p$. There exist arbitrarily small numbers $\rho > 0$ such that $\{|z-p| = \rho \}$ does not meet $\overline{U} \smallsetminus \overline{\mathbb D}$, hence for those $\rho$ the set $\{|z-p| = \rho \} \cap \{|z| > 1  \}$ is contained in the complement of $\overline{U} \smallsetminus \overline{\mathbb D}$ in $\mathbb C \smallsetminus \overline{\mathbb D}$, namely in Int$V \cap  \{|z| > 1  \}$. There cannot be two disjoint open connected subsets of $\mathbb C \smallsetminus \overline{\mathbb D}$ for which $p$ is an accumulation point, which both contain half-circles $\{|z-p| = \rho \} \cap \{|z| > 1  \}$ for some positive numbers $\rho$. Since the open set $\mathbb C \smallsetminus ( \overline{\mathbb D} \cup \overline{U})$ is not thin at $p$ it has exactly one such component and this component is non-thin at $p$. Theorem \ref{Thm1} is proved.
\end{proof} 
The proof of Theorem \ref{Thm1a} is a slight modification of the proof of Theorem \ref{Thm1}. We will omit it.
\begin{Prop}\label{Propp3}
Suppose the continuous function $\mathcal{F}$ on a closed neighborhood $V=D(p,r) \smallsetminus U$ of $p$ is finely analytic at $p$. Suppose $\gamma : [-1,1] \rightarrow \mathbb C$ is a smooth arc with $\gamma (0) = p$ which divides $\overline{D(p,r)}$ into two connected components $D_+(p,r)$ and $D_-(p,r)$. Then there is a smaller fine neighborhood $V_1$ of $p$ such that $\mathcal{F}|_{D_+(p,r) \cap V_1}$ is uniquely determined by $\mathcal{F}|_{D_-(p,r) \cap V}$.
\end{Prop}
\begin{proof}
It is enough to show that if $\mathcal{F}$ is finely analytic at $p$ on $V$ and $\mathcal{F}|_{D_-(p,r) \cap V} \equiv 0$  then $\mathcal{F}|_{V_1}$ is equal to zero for some fine neighborhood $V_1 \subset V$ of $p$. As in the proof of Theorem \ref{Thm1} there exist compact subsets $K_n$ of $U$, each being the finite disjoint union of closures of simply connected domains and analytic functions $\mathcal{F}_n$ on $D(p,r) \smallsetminus K_n$ which are continuous on $\overline{D(p,r) \smallsetminus K_n}$ and uniformly bounded by a constant $C>1$, and converge to $\mathcal{F}$ uniformly on $V$. 
Let $J$ be a closed arc of a circle $\partial D(p,\rho)$ for some $\rho >0$ which is contained in $D_-(p,r) \cap V$ and has length at least $5 \pi \rho /6$. Since $\mathcal{F}=0$ on $J$, the numbers $\epsilon_n \overset{def}{=} \max_J|\mathcal{F}_n|$ are less than $1$ for $n > n_0$ and tend to zero for $n \rightarrow \infty$.
The same arguments as in the proof of Theorem \ref{Thm1} give a number $r_1 >0$ and an open set $U_1$ which is thin at $p$ and a harmonic measure estimate analogously to (\ref{hm1}) such that the Two-constant Theorem gives for $z \in V_1 = \overline{D(p,r_1)} \smallsetminus U_1$ and all $n > n_0$
\begin{eqnarray*}
\log |\mathcal{F}_n(z)| & \leq & \log  \epsilon_n \cdot \omega(z,J, D(p,\rho) \smallsetminus K_n) + \\
&+& \log C \cdot (1- \omega(z,J, D(p,\rho) \smallsetminus K_n) ) \leq \\
&\leq&  \log  \epsilon_n \cdot \frac{1}{4} +  \log C. 
\end{eqnarray*}
Hence $\mathcal{F}(z) = \lim_{n \rightarrow \infty} \mathcal{F}_n(z) = 0$ for $z \in V_1$.
\end{proof}

\begin{proof}[Proof of Corollary \ref{Cor1}]
By Proposition \ref{Prop1}, $f|_{D_i}$ has fine analytic continuation to a set $V_1 = D(p,r_1) \smallsetminus U_1 \subset \overline{D}_i \cup \overline{D}_e$ for some $r_1 > 0$ where $U_1$ is an open set which is thin at $p$. Moreover, this fine analytic continuation equals $f|_{V_1}$. The domain $D_e$ is non-thin at $p$, since $D_e$ contains $D(p, \rho) \cap \{|z| > 1  \} \smallsetminus \overline{U}$, where $\overline{U}$ is thin at $p$. By Theorem \ref{Thm1} $\Gamma_f(V_1)$ is in the pluripolar hull of $\Gamma_f(D_i)$, and hence, as in the proof of the second part of Theorem \ref{Thm1}, $\Gamma_f(D_e)$ is contained in the pluripolar hull of $\Gamma_f(D_i)$.
\end{proof}
\begin{proof}[Proof of Corollary \ref{Cor2}:] Let $\bigcup_l I_l$ be a union of disjoint open
arcs on $\partial \mathbb D$ such that $\bigcup_l I_l$ is dense on
$\partial \mathbb D$, its linear measure is $< 2 \pi$ and  $\bigcup_l
I_l$ is thin at $1$. Let $G \subset \mathbb D$  be a Jordan domain 
 whose boundary $\gamma$ is a smooth, nowhere
analytic curve with  $\partial \mathbb D \smallsetminus \gamma 
\subset \bigcup_l \lozenge_l$, where the $\lozenge_l$:s are similar rhombs corresponding to the arcs $I_l$ like in Section $\ref{intro}$.
Let $f$ be a conformal mapping of $G$ onto $\mathbb D$. $f$ extends to a
smooth homeomorphism of $\overline{G}$ onto $\overline{\mathbb D}$. 
Since $\gamma$ is nowhere analytic
$f$ and its inverse $f^{-1}$ do not have analytic continuation across 
any part of the boundary of their domain of definition. However by Schwarz reflection principle $f$ admits
pseudocontinuation across the set $\mathcal{E} = \partial \mathbb D \smallsetminus \bigcup_l I_l$ and hence extends to a function  in $A(\overline{D_i} \cup \overline{D_e})$ for suitable domains
$D_i$ and  $D_e$ of the kind described before the statement of Theorem $\ref{Thm1}$. Note that the extended function is also univalent. By
Corollary \ref{Cor1} the
pluripolar hull of $\Gamma_{f}(D_i)$ (hence of $\Gamma_{f}(G)$) contains 
$\Gamma_{f}(D_e)$. The graph of $f$ over the subset $G$ in the $z$-plane,
$\Gamma_f(G)= \{ (z,w) \in \mathbb C^ 2 \vert z \in G, w=f(z) \}$, can be
considered as the graph of its inverse function over the set
$\mathbb D$ in the $w$-plane, $ \{ (z,w) \in \mathbb C^2 \vert w \in \mathbb D , z=f^{-1}(w) \}$.
The Corollary follows.
\end{proof}

\section{Points which are not in the pluripolar hull}\label{Section3}
In this Section we will prove Theorem \ref{Thm2}. For the proof it will be convenient to use the following known results.
 
Let $\Omega$ be a pseudoconvex domain in $\mathbb C^N$. In \cite{hull} the {\em{negative pluripolar hull }} is defined as
\begin{displaymath}
E_{\Omega}^- \overset{def}{=} \bigcap \{z \in \Omega : u(z)= -\infty \},
\end{displaymath}
where the intersection is taken over all{\em{ negative
    }} plurisubharmonic functions in $\Omega$ that are $-\infty$ on
$E$. 
 The following relation between the negative pluripolar hull and the
 pluripolar hull holds (see \cite{hull}).
\begin{Thm}\label{ww}
Let $\Omega$ be a pseudoconvex domain in $\mathbb{C}^N$. Let $\{\Omega_j\}$ be an increasing sequence of relatively compact subdomains of $\Omega$ with $\bigcup_j \Omega_j = \Omega$. Let $E \subset \Omega$ be pluripolar. Then
\begin{displaymath}
E_{\Omega}^* = \bigcup_j(E \cap {\Omega}_j)_{\Omega_j}^-.
\end{displaymath}
\end{Thm}

For a subset $E \subset \Omega$, the {\em{pluriharmonic measure}} at a point $z \in \Omega$ of $E$ relative to $\Omega$, is defined as 
\begin{equation}\label{measure}
 W(z,E,\Omega)= -\sup\{u(z) : u  \textrm{ is plurisubharmonic in } \Omega \textrm{ and u} \leq -\chi_E \},
\end{equation}
where $\chi_E$ is the characteristic function of the set $E$. The relation between the negative pluripolar hull and pluriharmonic measure is given in the following Theorem \cite{hull}.
\begin{Thm}\label{Thm4}
Let $\Omega$ be a domain in $\mathbb{C}^N$ and let $E \subset \Omega$ be pluripolar. Then
\begin{displaymath}
E_{\Omega}^- = \{z \in \Omega : W(z,E,\Omega) > 0 \}.
\end{displaymath}
\end{Thm}
The following theorem was recently proved by Wiegerinck and Edigarian in \cite{Ediwie3}.
\begin{Thm}\label{Thm5}
Let $\Omega$ be a pseudoconvex open set in $\mathbb C^N$ and let $E \subset \Omega$ be a $F_\sigma$ pluripolar subset. Assume that $E$ is connected. Then $E^*_\Omega$ is also connected.
\end{Thm}
In \cite{zeta} A. Zeriahi proved the following Theorem:
\begin{Thm}[Zeriahi]
Let $\Omega$ be a pseudoconvex domain and $F$ a pluripolar subset of $\Omega$ of type $F_\sigma$ and let $E$ be a closed subset of $\Omega$ such that $E \supset F_\Omega^*$. Then there exists a plurisubharmonic function $u$ on $\Omega$ which is continuous and $> -\infty$ on $\Omega \smallsetminus E$ and $= -\infty$ on $F$.
\end{Thm}
For a natural number $j$ we denote by $\mathbb B_j$ the (open) ball of radius $j$ and center $0$ in $\mathbb C^2$.
The key in the proof of Theorem \ref{Thm2} is contained in the following Lemma.
\begin{Lem}\label{Lem2}
Let $f$ be analytic in $\mathbb D$ and let $K$ be a closed disk in $\mathbb D$. Choose $j_0$ so that $\Gamma_f(K) \subset \mathbb B_{j_0}$.  Put $E = \overline{(\Gamma_f(\mathbb D))^*_{\mathbb C^2}}$, let $j \geq j_0$ and consider the compact set $E_j = E \cap \overline{\mathbb B_j}$ and the open set $U_j = \mathbb C \smallsetminus \pi_1(E_j)$. Suppose $p \in \partial \mathbb D$ is non-thin for $U_j$. Then for each $\epsilon > 0$ there exists a continuous plurisubharmonic function $g$ in $\mathbb B_j$, such that $g \leq 0$, $g=-1$ on $\Gamma_f(K)$ and $g(p,w) \geq - \epsilon$ for any $w \in \mathbb C$ such that $(p,w) \in \mathbb B_j$. In particular
\begin{displaymath}
W((p,w), \Gamma_f(K), \mathbb B_j ) = 0
\end{displaymath}
for each $w \in \mathbb C$ such that $(p,w) \in \mathbb B_j$.
\end{Lem}
\begin{Cor}\label{Cor5}
Let $f$ be analytic in $\mathbb D$ and let $p \in \partial \mathbb D$. Then either $\{p \} \times \mathbb C$ does not meet $(\Gamma_f(\mathbb D))^*_{\mathbb C^2}$ or $\pi_1($ $\overline{(\Gamma_f(\mathbb D))^*_{\mathbb C^2}}$ $)$ is a fine neighborhood of $p$.
\end{Cor}
The Corollary improves a result of \cite{Ediwie3}.
\begin{proof}[Proof of Corollary \ref{Cor5}]
Suppose $\pi_1(E) = \pi_1($ $\overline{(\Gamma_f(\mathbb D))^*_{\mathbb C^2}}$ $)$ is not a fine neighborhood of $p$, i.e. $U= \mathbb C \smallsetminus \pi_1(E)$ is not thin at $p$. Let $K$ be a closed disk contained in $\mathbb D$ and let $j_0$ be so large that $\Gamma_f(K) \subset \mathbb B_{j_0}$. Let $j \geq j_0$. Since $E_j= E \cap \overline{\mathbb B_j} \subset E$ the set $\pi_1(E_j)$ is also not a fine neighborhood of $p$, i.e. the set $U_j = \mathbb C \smallsetminus \pi_1(E_j)$ is not thin at $p$. By Lemma \ref{Lem2} and Theorem \ref{Thm4} for each $j \geq j_0$ the set $ (\{p\} \times \mathbb C) \cap \mathbb B_j$ does not meet the set $(\Gamma_f(K))^-_{\mathbb B_j}$. By Theorem \ref{ww} $\{p\} \times \mathbb C$ does not meet $(\Gamma_f(K))^*_{\mathbb C^2}$. Since the $\mathbb C^2$-pluripolar hull of $\Gamma_f(K)$ and $\Gamma_f(\mathbb D)$ coincide the Corollary is proved.
\end{proof}

Note that the arguments of Lemma \ref{Lem2} may be applied in situations when $\mathbb D$ is replaced by another planar domain. In particular if $\pi_1(E)$ is the closure of a domain, this gives a tool to study $(\Gamma_f(\mathbb D))^*_{\mathbb C^2}$ over boundary points of the domain. We will not work this out here.
\begin{proof}[Proof of Theorem \ref{Thm2}]
Suppose, on the contrary, that for each point $p \in \partial \mathbb D$ the set $U= \mathbb C \smallsetminus \pi_1(E)$ is non-thin at $p$. By Corollary \ref{Cor5} the set $\{p\} \times \mathbb C$ does not meet $(\Gamma_f(\mathbb D))^*_{\mathbb C^2}$ for each $p \in \partial \mathbb D$. Hence by Theorem \ref{Thm5} $(\Gamma_f(\mathbb D))^*_{\mathbb C^2} \subset \mathbb D \times \mathbb C$. 
\end{proof}
\begin{proof}[Proof of Lemma \ref{Lem2}]
In case $p \in U_j$ we decrease $U_j$ by replacing it by a simply connected set $\widetilde{U}_j$ with $p \in \partial \widetilde{U}_j$ which is still not thin at $p$ (e.g. cut the connected component of $U_j$ containing $p$ along a curve joining $p$ with the boundary $\partial U_j$). Put $\widetilde{U}_j = U_j$ if $p \notin U_j$. Then $\widetilde{U}_j \subset U_j$ is not thin at $p$, $p \notin \widetilde{U}_j$ and $\widetilde{U}_j$ is simply connected.
Let $K_n \subset \widetilde{U}_j$ be an increasing sequence of compacts, each of them being the finite disjoint union of closures of smoothly bounded simply connected domains, such that $\widetilde{U}_j = \bigcup K_n$. Then for each natural number $n$ the set $D_n = \widehat{\mathbb C} \smallsetminus (K \cup K_n)$ is a domain. Note that $D_n$ has regular boundary for the Dirichlet problem. Hence the harmonic measure $\omega(z,\partial K, D_n)$ of $ \partial K$ for the domain $D_n$ computed at the point $z\in \overline{D}_n$ is a continuous function on $\overline{D}_n$ which is harmonic on $D_n$.

We claim that 
\begin{equation}\label{limit}
\lim_{ n \rightarrow \infty} \omega(p,\partial K, D_n) =0.
\end{equation}
Indeed, extend for each $n$ the function $\omega(z,\partial K, D_n)$ to the set $K_n$ by putting it equal to zero there. Denote the extended function by $h_n$. $h_n$ is continuous and subharmonic on $\widehat{\mathbb C} \smallsetminus K$, it is non-negative and $h_{n+1} \leq h_n$ for each $n$. Hence $h = \lim_{n \rightarrow \infty} h_n$ is non-negative and subharmonic on $\widehat{\mathbb C} \smallsetminus K$ (being the decreasing limit of a sequence of subharmonic functions). Moreover, $h=0$ on $\widetilde{U}_j$. Since $\widetilde{U}_j$ is non-thin at $p$, $p$ is an accumulation point of $\widetilde{U}_j$ and $h(p)=0$. Hence the claim.
 
Let now $\epsilon > 0$ and choose $n$ so that 
\begin{equation}\label{hm}
\omega(z,\partial K, D_n) < \epsilon.
\end{equation}
 
The function
\begin{equation}
v(z,w) = \left \{ \begin{array}{ll}
        -\omega(z,\partial K,D_n), & z \in D_n, (z,w) \in \mathbb B_j \\
        -1, & z \in K, (z,w) \in \mathbb B_j \\
        \end{array} \right.
\end{equation}
is plurisubharmonic on a large part of the ball $ \mathbb B_j$, precisely on $\{(z,w) \in \mathbb B_j: z \notin K_n  \}$. We want to obtain a plurisubharmonic function in the whole $\mathbb B_j$ using a standard gluing procedure near $(\partial K_n \times \mathbb C) \cap \mathbb B_j$ together with the Theorem of Zeriahi. Denote by
$u$ a plurisubharmonic function on $\mathbb B_{j+1}$ which is
continuous and $> -\infty$ on $\mathbb B_{j+1} \smallsetminus E$ (with
$E= \overline{(\Gamma_f(\mathbb D))^*_{\mathbb C^2}}$ ) and equals $- \infty$ on $\Gamma_f(\mathbb D)$. Its existence is guaranteed by Zeriahi's Theorem. Adding a constant we may assume that $u<0$ on the compact set $\overline{\mathbb B_j} \subset \mathbb B_{j+1}$. Let $G_n \subset \widetilde{U}_j$ be a smoothly bounded open set such that $K_n \subset G_n \subset \overline{G_n} \subset \widetilde{U}_j$. (We think of $G_n$ being close to $K_n$.) Let $V_n$ be a small neighborhood of $\partial G_n$ such that $\partial G_n \subset V_n \subset \overline{V_n}\subset \widetilde{U}_j \smallsetminus K_n \subset D_n$. Then there exists a positive constant $\delta$ such that $- \omega(z,\partial K, D_n) < - \delta$ on $\overline{V_n}$. Since $\overline{V_n} \subset U_j = \mathbb C \smallsetminus \pi_1(E_j)$, the set $(\overline{V_n} \times \mathbb C) \cap \overline{\mathbb B}_j$ does not meet $E$, hence the function $u$ is $> - \infty$ and continuous on this set. Multiplying $u$ by a positive constant we may assume that $u> -\delta$ on the compact subset  $(\overline{V_n} \times \mathbb C) \cap \overline{\mathbb B_j}$. 
For $(z,w) \in \mathbb B_j$ define
\begin{equation}
g(z,w) = \left \{ \begin{array}{ll}
        u(z,w), & z \in G_n, (z,w) \in \mathbb B_j \\
        \max \{v(z,w),u(z,w) \}, & z \notin G_n, (z,w) \in \mathbb B_j. \\
        \end{array} \right.
\end{equation}
$v$ is defined on  $\{(z,w) \in \mathbb B_j: z \notin K_n \}$ and $K_n \subset G_n$, hence $g$ is well defined. On $(\overline{V_n} \times \mathbb C) \cap \mathbb B_j$ we have the inequality $v < -\delta < u$, hence $g=u$ on $(V_n \times \mathbb C)  \cap \mathbb B_j$. Since $u$ and $v$ are plurisubharmonic where they are defined, the function $g$ is plurisubharmonic on $\mathbb B_j$. Since for $(z,w) \in \Gamma_f(K) \cap \mathbb B_j$ the relations $u(z,w) = - \infty$, $v(z,w) = -1$ hold we obtain for these points $g(z,w) =-1$.
On the other hand, since $p \notin G_n$, for points of the form $(p,w) \in \mathbb B_j$ we have by (\ref{hm})
\begin{displaymath}
g(p,w) \geq v(p,w) > - \epsilon
\end{displaymath}
The Lemma is proved.
\end{proof}

\end{document}